\documentclass[11pt]{amsart}
\usepackage{latexsym}
\usepackage{amsmath}

\title[Algebraically Constructible Functions]
{Algebraically Constructible Functions\\
and signs of polynomials}
\author{Adam Parusi\'nski}
\address{D\'epartement de Math\'ematiques, Universit\'e d'Angers, 2 bd.
Lavoisier, 49045 Angers Cedex, France, and School of Mathematics and
Statistics, University of Sydney, Sydney, NSW 2006, Australia}
\email{parus@tonton.univ-angers.fr, parusinski\_a@maths.su.oz.au}
\author{Zbigniew Szafraniec}
\address{Institute of Mathematics, University of Gda\'nsk,
Wita Stwosza 57, 80-952 Gda\'nsk, Poland}
\email{szafran@ksinet.univ.gda.pl}
\subjclass   {Primary: 14P05, 14P25 }

\keywords {real algebraic set, Euler characteristic, algebraically 
constructible functions}

\thanks {Research partially supported by an Australian Research Council Small 
Grant.  Second author also partially supported by KBN 610/P3/94. }

\newcommand{\abstracttext}{}

\newcommand{\calI}{{\mathcal I}}
\newcommand{\calS}{{\mathcal S}}

\newfont{\bboard}{msbm10 scaled 1100}

\renewcommand{\labelenumi}{(\roman{enumi})}

\newcommand{\prawo}{\mbox{$\longrightarrow$}}

\newcommand{\eR}{\mbox{$\mathbf R$}}

\newcommand{\Apis}{\mbox{$\mathcal A$}}

\newcommand{\Ipis}{\mbox{$\mathcal I$}}
\newcommand{\Npis}{\mbox{$\mathcal N$}}

\newcommand{\Rn}{\mbox{${\mathbf R}^{n}$}}  
\newcommand{\Rm}{\mbox{${\mathbf R}^{m}$}}
\newcommand{\Cn}{\mbox{${\mathbf C}^{n}$}}

\newcommand{\RndoR}{\mbox{${\mathbf R}^{n}$}
\mbox{$\longrightarrow$}\mbox{$\mathbf R$}}

\newcommand{\RndoRn}{\mbox{$ {\bf R}^{n} $}
\mbox{$ \longrightarrow $}\mbox{$ {\mathbf R}^{n} $}}

\newcommand{\gO}{\mbox{${\mathbf 0}$}}
\newcommand{\go}{\mbox{$\mathbf 0$}}

\newcommand{\eNen}{\mbox{\( {\mathbf N}^{n} \)}}

\newcommand{\ckrop}{\mbox{$\,\cdots\,$}}
\newcommand{\pkrop}{\mbox{$\,,\ldots,\,$}}
\newcommand{\plkrop}{\mbox{$\,+\cdots+\,$}}

\newcommand{\sign}{\mbox{$\operatorname {sgn}$}\,}
\newcommand{\sgn}{\mbox{$\operatorname {sgn}$}\,}
\newcommand{\supp}{\mbox{$\operatorname {supp}$}\,}
\newcommand{\signature}{\mbox{$\operatorname {signature}$}\,}
\newcommand{\one}{\mbox{$\mathbf 1$}}

    \newtheorem{theorem}{Theorem}[section]
    \newtheorem{lemma}[theorem]{Lemma}
    \newtheorem{proposition}[theorem]{Proposition}
    \newtheorem{corollary}[theorem]{Corollary}
 
   \vspace{4in}

\begin{document}
\begin{abstract} \abstracttext 
Let $W$ be a real algebraic set.  We show that the following families of 
integer-valued functions on $W$ coincide: (i) the functions of the form 
$w\to \chi(X_w)$, where $X_w$ are the fibres of a regular morphism 
$f:X\to W$ of real algebraic sets, (ii) the functions of the form  
$w\to \chi(X_w)$, where $X_w$ are the fibres of a proper regular morphism 
$f:X\to W$ of real algebraic sets, 
(iii) the finite sums of signs of polynomials on $W$.  Such 
functions are called algebraically constructible on $W$.  
Using their characterization in terms of signs of polynomials we present new 
proofs of their basic functorial properties with respect to the link operator 
and specialization.
\end{abstract} \maketitle
%\tableofcontents
%%%% **** The text of the paper starts here **** %%%% %
%
%       Intro
%

\section{Introduction}\label{introduction}
Let $f:X\to W$ be a regular morphism of real algebraic sets.  
Consider on $W$ an integer-valued function $\varphi(w) = \chi (X_w)$, which 
associates to $w\in W$ the Euler characateristic of the fibre $X_w=f^{-1}(w)$.  
The main purpose of this paper is to study the properties of such $\varphi$.  

Firstly, by stratification theory, $\varphi$ is (semialgebraically) 
constructible, that is there exists a semialgebraic stratification $\calS$ 
of $W$ such that $\varphi$ is constant on strata of $\calS$.  Equivalently, 
we may express this property by saying that $\varphi$ is bounded and 
$\varphi^{-1}(n)$ is semialgebraic for every integer $n$.  However it is 
well-known that not all semialgebraically constructible functions  on $W$ 
are of the form $\chi (X_w)$ for a regular morphism 
$f:X\to W$. For instance, if 
$W$ is irreducible, then $\chi (X_w)$  
has to be generically 
constant modulo 2, see for instance \cite[Proposition 2.3.2]{akbulutking}.  
Also in the case of $W$ irreducible, as shown in \cite{costekurdyka2},   
there exists a real polynomial $g:W\to \eR$ such that generically on $W$
$\varphi (w) \equiv \sgn g(w) \pmod 4$, where by $\sgn g$ we
denote the sign of $g$.   
As we show in Theorem \ref{wtyczka} below, for any regular morphism 
$f:X\to W$ of real algebraic sets there exist real polynomials
$g_1,\ldots, g_s$ on $W$ such that for every $w\in W$
\[ \chi (X_w)= \sgn g_1(w) + \cdots +  \sgn g_s(w) .\]  
In particular, taking  $g=g_1 \cdot  sg_s$ we recover
 the result of \cite{costekurdyka2}.   

Constructible  functions of the form $\varphi (w) = \chi (X_w)$, for 
proper regular morphisms $f:X\to W$, were studied in 
\cite{mccroryparusinski2} in a different context.  Following 
\cite{mccroryparusinski2}  
we call them {\it algebraically constructible}.  As shown in 
\cite{mccroryparusinski2} the 
family of algebraically constructible functions is preserved by various 
natural geometric operations such as, for instance, push-forward, duality, 
and specialization.  In a way they behave similarly to constructible 
functions on complex algebraic varieties.  However, unlike their complex 
counterparts, they cannot be defined neither in terms of stratifications 
nor as combinations of characteristic functions of real algebraic 
varieties, cf. \cite {mccroryparusinski2}.  Algebraically constructible 
functions  were used in 
\cite {mccroryparusinski2} to study local topological properties 
of real algebraic sets.  

In particular, Theorem \ref{wtyczka} below can be reformulated as follows.  
Algebraically constructible functions on a real algebraic set $W$ 
coincide with finite sums of signs of real polynomials on $W$, see Theorem 
\ref{key} below.  Using this characterization, in section 6, we give new, 
alternative proofs of basic properties of 
algebraically constructible functions, without using the resolution 
of singularities as in \cite{mccroryparusinski2}.    

The main result of the paper, Theorem \ref{wtyczka}, is proven in section 
\ref{families}.   In sections 2-4 we develop necessary techniques for 
the proof and recall basic results the proof is based on.    
In particular, in section \ref{preliminaries} we recall   
the Eisenbud-Levine Theorem \ref{el} and the Khimshiashvili formula 
\ref{deszczyk}.  In section \ref{division} we review some basic facts on 
the Grauert-Hironaka formal division algorithm Theorem \ref{mrowka}, which
we then use to obtain a parametrized version of the Eisenbud-Levine 
Theorem, Propositions \ref{trasa} and \ref{tara}, with parameter in a given 
algebraic set $w\in W$.  
In section \ref{vectorfields} we study polynomial families 
of polynomial vector fields  $F_w:\Rn\prawo\Rn$ parametrized by $w\in W$.  

The proof of Theorem \ref{wtyczka} can be sketched briefly as follows.  
First, by an argument similar to the Khimshiashvili Formula, we show that the 
Euler characteristic $\chi(X)$ of a real algebraic set $X$ can be calculated 
in terms of the local topological degree at the origin of a polynomial 
vector field, see Proposition \ref{guzik}.  Then using the theory developed 
by the second named author, see e.g. \cite{szafraniec4,szafraniec17},
  we generalize this observation 
in two directions.  Firstly, we show that 
for a regular morphism $f:X\to W$,  
the Euler characteristic of the fibers $\chi(X_w)$ can be 
expressed in terms of the local topological degree $\deg_{0} G_w$ 
at the origin of a 
family $G_w:\Rn\prawo\Rn$ of polynomial vector fields, which depends 
polynomially on $w$.  Secondly, 
as shown in Lemma \ref{marazm}, we may choose all $G_w$ in 
such a way that they have algebraically isolated zero at the origin.  Then, 
by the Eisenbud-Levine Theorem \ref{el}, each $\deg_{0} G_w$ can be calculated 
algebraically, that is $\deg_{0} G_w$ equals the signature of an associated 
symmetric bilinear form $\Psi_w$.  By section 3, we may as well require that 
$\Psi_w$ depend "polynomially" on $w$.  More precisely, there exists a 
symmetricmatrix $T(w)$ (representing $\Psi_w$) 
with entries polynomials in $w$, such that $\deg_{0} G_w$ equals 
the signature of $T(w)$, for all $w$ in a Zariski open subset of $W$.   
(See \ref{trasa} and \ref{tara} for the details.)   
Finally by Descartes' Lemma, we express the signature 
of $T(w)$  in terms of signs of polynomials in $w$, 
see Lemma \ref{rosa} and the proof of Lemma 
\ref{sloniczek},.  

For the definitions and properties of real algebraic 
sets and maps we refer the reader to
\cite{benedettirisler}. 
By a real algebraic set we mean the locus of zeros of a finite set of
polynomial functions on $\bf R^n$.

\medskip
\section{Preliminaries}\label{preliminaries}
Let $f(x)=a_n x^n+a_{n-1}x^{n-1}\plkrop a_0$ be a real polynomial.
Let $\Lambda$ be the set of all pairs $(r,s)$
with $0\leq r<s\leq n$ such that $a_r\neq 0, a_s \neq 0$, and
$a_i=0$ for $r<i<s$.  Denote $\Lambda'=\{(r,s)\in\Lambda\mid r+s
\mbox{ is odd }\}$.

\begin{lemma}\label{rosa}
Assume that all roots of $f(x)$ are real and $a_0\neq 0, a_n\neq 0$.  
Let $p_+$ (resp. $p_-$) denote the number of positive (resp. negative)
roots counted with multiplicities. Then
\[p_+-p_-=-\sum {\sgn}\, a_r a_s,\mbox{ where }(r,s)\in \Lambda',\]
\[p_+-p_-\equiv n+1+(-1)^{n+1}{\sgn}\, a_0a_n\pmod{4}.\]
\end{lemma}

{\em Proof.\/} We say that the pair of real numbers $(a, b)$ changes
 sign if $ab<0$. If this is the case then
$(1-\sign ab)/2=1$, if $ab>0$ then
$(1-\sign ab)/2=0$.

As a consequence of Descartes' lemma (see 
\cite [Theorem 6, p.232]{mostowskistark}, 
or \cite [Exercise 1.1.13 (4), p.16]{benedettirisler}), $p_+$ equals
the number of sign changes in the sequence of non-zero coefficients
of $f(x)$, that is 
\[ p_+=\sum(1-\sign a_r a_s)/2,\mbox{ where } (r,s)\in\Lambda.\]
According to the same fact, $p_-$ equals  the number of sign changes 
in the sequence of non-zero coefficients of $f(-x)$, i.e.
\[ p_-=\sum (1-(-1)^{r+s}\sign a_r a_s)/2, \mbox{ where }(r,s)\in\Lambda.\]
Hence
\[p_+-p_-=-\sum \sign a_r a_s, \mbox{ where } (r,s)\in\Lambda'.\]

The sign of the product of all roots, that is $(-1)^{p_-}$, equals
$(-1)^n\sign a_0 a_n$. Thus
$2p_-\equiv 3+(-1)^{p_-}=3+(-1)^n\sign a_0 a_n\pmod{4}$.
Finally, since $p_++p_-=n$, we conclude that 
\[p_+-p_-=n-2p_-\equiv n-3-(-1)^n\sign a_0 a_n
   \equiv n+1+(-1)^{n+1}\sign a_0 a_n\pmod{4}. \Box  \] 
\smallskip

 % This is file: deform.tex
Let $F:(\Rm,\go)\prawo(\Rm,\go)$ be a germ of a continuous mapping
with  isolated zero at $\go$.   Then we denote by $\deg_0 F$ 
 the local topological degree of $F$ at the origin.
Suppose, in addition, that $F=(f_1\pkrop f_m)$ is a real analytic germ. Let
$\eR [[x]]=\eR [[x_1\pkrop x_n]]$ denote the ring of formal power
series and let $I$ denote the ideal in $\eR[[x]]$ generated by 
$f_1\pkrop f_m$.  Then $Q=\eR[[x]]/I$ is an $\eR$-algebra.  
If $\dim_{\bf R} Q<\infty$,  
then $\go$ is isolated in $F^{-1}(\go)$ and in this case we say that
$F$ has {\em an algebraically isolated zero at\/} $\go$.  Let $J$
denote the residue class in $Q$ of the Jacobian determinant 
\[\frac{\partial(f_1\pkrop f_m)}{\partial(x_1\pkrop x_m)} . \]
The next theorem is due to Eisenbud and Levine \cite{eisenbudlevine}, 
see also \cite{arnoldetal,beckeretal,khimshiashvili} for a proof.  
\smallskip

\begin{theorem}[Eisenbud\&Levine Theorem]
\label{el} Assume that $\dim_{ \bf R} Q<\infty$. Then
\begin{enumerate}
\item $J\neq 0$ in $Q$,
\item for any $\eR$-linear form $\varphi:Q\prawo\eR$ such that $\varphi(J)>0$, 
the corresponding symmetric bilinear form 
$\Phi: Q\times Q\prawo\eR$, $\Phi(f,g)=\varphi(fg)$, is 
non-degenerate and
\[\signature \,\Phi=\deg_0 F.\quad \Box\]
\end{enumerate}\end{theorem}

The next formula was proved by Khimshiashvili \cite{khimshiashvili},
for other proofs see \cite{arnold,arnoldetal,wall}.

\begin{theorem}[Khimshiashvili Formula] \label{deszczyk}
Let $g:(\Rm,\gO)\prawo(\eR,0)$ be a real analytic germ  with isolated
critical point at $\go$.  Let $S_\epsilon$ denote the sphere of a small radius
$\epsilon$ centered at the origin and let $A_\epsilon=S_\epsilon
\cap\{g\leq 0\}$. (Note that all $A_\epsilon$ are
homeomorphic for $\epsilon>0$ small enough.)  
Then the gradient $\nabla g:\Rm\prawo\Rm$ of $g$ has an
isolated zero at $\go$ and
\[\chi(A_\epsilon)=1-\deg_0(\nabla g).\hspace{1.0cm}\Box\]
\end{theorem} 
\medskip

\begin{lemma}  \label{biurko}
Let $g:\Rn\times\eR\prawo \eR$ be a polynomial vanishing at $\go$ 
and such that if $g(x,t)\leq0$ then either
$(x,t)=\gO$ or $t>0$.
Let $S_\epsilon\subset\Rn\times\eR$ (resp. $B_\epsilon$)
denote the sphere (resp. the open ball) of radius $\epsilon$
centered at the origin and let
$A_\epsilon=S_\epsilon\cap\{g\leq 0\}$. Let
$P_\eta=\Rn\times\{\eta\}$ and
$M_{\epsilon,\eta}=P_\eta\cap\{g\leq 0\}\cap B_\epsilon$.
Then, for $0<\eta \ll\epsilon\ll 1$, $A_\epsilon$ and
$M_{\epsilon,\eta}$ have the same homotopy type. In particular,
\[\chi(A_\epsilon)=\chi(M_{\epsilon,\eta}).\hspace{1.0cm} %\\[1em]
\] 
\end{lemma}

{\em Proof.\/}  Consider on $N=\{(x,t)\mid g(x,t)\leq0\}$,  
the functions $\omega_1(x,t)=\|  x\| ^2+t^2$ and
$\omega_2(x,t)=t$.  Both $\omega_1$ and $\omega_2$ are
non-negative on $N$ and
$\omega_{1}^{-1}(0)\cap N=\omega_{2}^{-1}(0)\cap N=\{\go\}$.  
Let $N_{i}^{y}=\{(x,t)\in N\mid 0<\omega_i(x,t)\leq y\}$.  
By the topological local triviality of semi-algebraic mappings,  
see for instance \cite[Theorem 9.3.1]{bochnaketal}  or
\cite{coste1,hardt},  there is $\delta>0$ such that
$\omega_i:N_{i}^{\delta}\prawo (0,\delta]$, $i=1,2$, are topologically
 trivial fibrations.  For $0<y\leq \delta$ let $M_{2}^{y}$ denote the union
of connected components of $N_{2}^{y}$ containing $\go$
in their closures. Then
$\omega_2:M_{2}^{\delta}\prawo (0,\delta]$ is also 
topologically trivial.

Hence there exist constants $0<\alpha<\beta<\gamma<\delta$ such that
$M_{2}^{\alpha}\subset N_{1}^{\beta}\subset M_{2}^{\gamma}
\subset N_{1}^{\delta}$.   By the topological triviality,
the inclusions $M_{2}^{\alpha}\subset M_{2}^{\gamma}$ and
$N_{1}^{\beta}\subset N_{1}^{\delta}$ are homotopy equivalencies 
and hence so are $M_{2}^{\alpha}\subset N_{1}^{\beta}$
and $N_{1}^{\beta}\subset M_{2}^{\gamma}$.  

By the above, the total spaces of fibrations 
$\omega_1:N_{1}^{\delta}\prawo (0,\delta]$, 
$\omega_2:M_{2}^{\delta}\prawo (0,\delta]$  
are homotopy
equivalent to their fibers.  Consequently the fibers of both fibrations 
are homotopy equivalent. Now, to complete the proof,
 it is enough to observe that these fibers are of the form 
$A_\epsilon$ and $M_{\epsilon,\eta}$, where $0<\eta\ll\epsilon\ll 1$. $\Box$

\begin{proposition} \label{guzik}
Let $f:\Rn\times\eR\prawo\eR$ be a non-negative homogeneous polynomial
of degree $2d$ such that $f(x,0)=\|  x\| ^{2d}$.
Let $X=\{x\in\Rn\mid f(x,1)=0\}$ and define 
$g(x,t)=f(x,t)-t^{2d+1}$. Then $g$ has an isolated critical
point at the origin and
\[\chi(X)=1-\deg_0(\nabla g).%\\[1em]
\]
\end{proposition} 

{\em Proof.\/} Let
\[\Sigma=\left\{(x,t)\mid \frac{\partial f}{\partial x_1}=
  \cdots=\frac{\partial f}{\partial x_n}=0\right\},\]
$P_\eta=\Rn\times\{\eta\}$, and 
$\Sigma_\eta=\Sigma\cap P_\eta$.
Let $f_\eta$ (resp. $g_\eta$) denote the restriction of $f$
(resp. $g$) to $P_\eta$.  Then $\Sigma_\eta$ is the set
of critical points of both $f_\eta$ and $g_\eta$.
We have $\Sigma_0=\{\gO\}$.  Since the set of critical values of 
any polynomial is finite, so is each $f_\eta(\Sigma_\eta)$.  Moreover, 
since $f$ is non-negative homogeneous of degree $2d$ and $\Sigma$ is 
a homogeneous set, there is $D>0$ such that any $y\in f_\eta(\Sigma_\eta)$, 
if non-zero, satisfies $y >D\mid\eta\mid^{2d}$.

If $\eta<0$, then $g_\eta>0$ and 
$0\in\eR$ is a regular value of $g_\eta$.  Clearly $g_0$ has a single
critical point at the origin.

Consider $0<\eta\ll 1$.  Let $x\in\Sigma_\eta$.  If 
$f_\eta(x)>0$ then
\[g_\eta(x)=f_\eta(x)-\eta^{2d+1}>D\eta^{2d}-\eta^{2d+1}>0.\]
If $f_\eta(x)= 0$ then
$g_\eta(x)<0$. Thus $0\in\eR$ is a regular value for $g_\eta$.
Hence there is an open neighbourhood 
$U\subset\Rn\times\eR$ of the origin such that $0\in\eR$ is a regular
value of $g$ on $U-\{\gO\}$, i.e. $g$ has an isolated
critical point at the origin.

For $\eta$ fixed, $\lim f_\eta(x)=+\infty$ as $\| x\| \prawo +\infty$.  
Denote
$N_\eta=\{x\mid f_\eta(x)= 0\}$ and
$M_\eta=\{x\mid g_\eta(x)\leq 0\}=\{x\mid f_\eta(x) \leq\eta^{2d+1}\}$.
If $\eta<0$, then $M_\eta$ is empty and $M_0=\{\gO\}$.
If $\eta>0$, then $N_\eta\subset M_\eta$.  As we have shown above, for
$0<\eta\ll 1$ both $f_\eta$ and $g_\eta$ have no critical points
in $M_\eta-N_\eta$.  Hence $N_\eta$ is a deformation retract of
$M_\eta$ and,  in particular, $\chi(N_\eta)=\chi(M_\eta)$.

Suppose $0<\eta\ll\epsilon$.  Then, since $f_0=g_0= \|x\|^{2d}$, 
$M_\eta\subset B_\epsilon$, that is $M_\eta=M_{\epsilon,\eta}$ 
in the notation of Lemma \ref{biurko}.  
Moreover, let $A_\epsilon=S_\epsilon\cap\{g\leq 0\}$.  By Lemma \ref {biurko}, 
$\chi(A_\epsilon)=\chi(M_\eta)$, and hence, by the above
\[ \chi(A_\epsilon)=\chi(M_\eta)= \chi(N_\eta). \] 
Finally, by the Khimshiashvili formula 
\ref{deszczyk}, 
\[\chi(A_\epsilon)=1-\deg_0 (\nabla g),\]
and the lemma follows since 
$\chi(X)=\chi(N_1)=\chi(N_\eta)$, for $\eta>0$.  $\hspace{1.0cm}\Box$

\bigskip
\section{The formal division algorithm}\label{division}

\noindent In the first part of this section we review some basic facts on 
the Grauert-Hironaka formal division algorithm for formal power series 
with polynomial coefficients.  In exposition and notation we follow closely 
 \cite{bierstonemilman}.  Then we apply the Grauert-Hironaka algorithm 
to derive a parametrized version of the Eisenbud-Levine Theorem \ref{el}, 
with parameter in a given algebraic set $W$.

Let $A$ be an integral domain. Let $A[[y]]=A[[y_1,\ldots,y_n]]$ denote
the ring of formal power series in $n$ variables with coefficients in $A$.

If $\beta=(\beta^1,\ldots,\beta^n)\in\eNen$, put
$\mid \beta\mid=\beta^1+\cdots+\beta^n$.
We order the $(n+1)$-tuples $(\beta^1,\ldots,\beta^n,\mid\beta\mid)$
lexicographically from the right. This induces a total ordering of $\eNen$.

Let $f\in A[[y]]$.  Write $f=\sum_{\beta\in N^n}f_{\beta}y^{\beta}$, where
$f_{\beta}\in A$ and $y^{\beta}$ denotes
$y_{1}^{\beta^1}\cdots y_{n}^{\beta^n}$. Let
$\supp (f)=\{\beta\in\eNen\mid f_\beta \neq 0\}$
and let $\nu (f)$ denote the smallest element of $\supp (f)$. Let
in$(f)$ denote $f_{\nu(f)}y^{\nu(f)}$.

Let $I$ be an ideal in $A[[y]].$
We define the diagram of initial exponents $\Npis (I)$ as
$\{\nu(f)\mid f\in I\}$. Clearly, $\Npis (I)+\eNen=\Npis (I)$.
There is a smallest finite subset $V(I)$ of $\Npis (I)$ such that
$\Npis (I)=V(I)+\eNen$. We call the elements of $V(I)$ 
the vertices of $\Npis (I)$.

Let $\beta_1,\ldots,\beta_t\in V(I)$ be the vertices of $\Npis (I)$ and 
choose 
$g^1\pkrop g^t\in I$ so that
$\beta_i=\nu(g^i)$,  $i=1,\ldots,t$.  The $\beta_1,\ldots,\beta_t$ 
induce the following decomposition of $\eNen$: Set
$\Delta_0=\emptyset$ and then define
$\Delta_i=(\beta_i+\eNen)\setminus \Delta_0\cup\ldots\cup\Delta_{i-1},\
i=1,\ldots,t$. Put
$\Delta=\eNen\setminus \Delta_1\cup\ldots\cup\Delta_t=\eNen\setminus \Npis(I)$.

Let $\mbox{in}(g^i)=g_{\beta_i}^i y^{\beta_i}$.
Then $g_{\beta_i}^{i}\neq 0$.  
Let $A_0$ denote the field of fractions of $A$. 
We denote by $S$ the multiplicative subset
of $A$ generated by the $g_{\beta_i}^{i}$ and by $S^{-1}A$ the
corresponding localization of $A$; i.e. the subring of $A_0$
comprising the quotients with denominators in $S$. Then
$S^{-1}A[[y]]\subset A_0[[y]]$.

\medskip
\begin{theorem}[Grauert, Hironaka,  
\cite{arocaetal,bierstonemilman,briancon,grauert}]\label{mrowka}
For every $f\in S^{-1}A[[y]]$ there exist unique
$g_i\in S^{-1}A[[y]]$, $\ i=1,\ldots,t$ , and
$r\in S^{-1}A[[y]]$ such that
$\beta_i +\supp (g_i)\subset \Delta_i$, $\supp (r)\subset\Delta$, and
\[f=\sum_{i=1}^{t}g_i g^i+r. \hspace{1.0cm}\Box\]
\end{theorem}

%The above theorem immediately implies

\begin{corollary}\label{rower}
$\nu(f)\leq\nu(r)$.  In particular, if $\Delta$ is finite and
$\beta<\nu(f)$ for all $\beta\in\Delta$, then $r=0$
and $f$ belongs to the ideal 
in $S^{-1}A[[y]]$ generated by $g^1,\ldots,g^t.\Box$
\end{corollary}

Let $S^{-1}I[[y]]$ denote the  ideal in $S^{-1}A[[y]]$ generated by $I$.
Then $S^{-1}A[[y]]/S^{-1}I[[y]]$ is finitely generated if and 
only if $\Delta$ is finite.
If this is the case then $S^{-1}A[[y]]/S^{-1}I[[y]]$ is a free
$S^{-1}A$ module and we take the monomials $y^\beta$, $\beta\in\Delta$,
as a basis.

Let $W\subset \Rn$ be an irreducible real algebraic  set and let $\Apis$ 
denote the ring of real polynomial functions on $W$.  Each $w\in W$ defines 
an evaluation homomorphism $h\mapsto h(w)$ of $\Apis$ onto $\eR$.  
For $f=\sum_{\beta}f_{\beta}y^{\beta}\in\Apis[[y]]$ we write
$f(x;y)=\sum_{\beta}f_{\beta}(x)y^{\beta}$,
and $f(w;y)=\sum_{\beta}f_{\beta}(w)y^{\beta}$ when the coefficients 
are evaluated at $x=w$.

Let $f^1\pkrop f^s\in\Apis[[y]]$ and let $\Ipis$ denote the ideal in
$\Apis[[y]]$ generated by $f^1\pkrop f^s$.  Let $\Npis =\Npis 
 (\calI)= \{\nu(g)\mid g\in\Ipis\}$ denote the diagram of initial exponents
(here $A=\Apis$).  Given $w\in W$.  We denote by$I_w$ 
 the ideal in $\eR[[y]]$ generated by
$f^1(w;y)\pkrop f^s(w;y)$ and by $\Npis_w=\Npis (I_w)$ the diagram
of initial exponents of $I_w$ (so here $A=\eR$).

The next theorem was proved by Bierstone and Milman 
\cite{bierstonemilman}.

\begin{theorem} \label{fotel}
Assume that $W$ is irreducible (so that $\Apis$ is an integral domain).
Let $\beta_1\pkrop\beta_t$ denote the vertices of $\Npis$ and 
choose $g^i\in\Ipis$ such that $\nu(g^i)=\beta_i$.  Let
\[\Sigma=\bigcup_{i=1}^{t}\{w\in W\mid g_{\beta_i}^{i}(w)=0\}.\]
Then $\Sigma$ is a proper algebraic subset of $W$, 
$\Npis_w=\Npis$ for all $w\in W-\Sigma$, 
$\nu(g^i)=\beta_i=\nu(g^i(w;\, \cdot\, ))$   for
every vertex $\beta_i\in \Npis$ and $w\in W \setminus \Sigma$. $\Box$
\end{theorem}

\begin{corollary}
Suppose that $\Delta_w=\eNen\setminus \Npis_w$ is finite for each 
$w\in W\setminus \Sigma$.
Then $\Delta=\eNen\setminus \Npis$ is also finite and $\Delta=\Delta_w$
for all $w\in W\setminus \Sigma$. $\Box$
\end{corollary}

Suppose that $\Delta$ is finite and let $\bar\beta$ denote 
the largest element in $\Delta$.  Let $j=y^{\bar\beta}$.
Then for $w\in W-\Sigma$, the residue class of $j$ in
$Q_w=\eR[[y]]/I_w$ is nonzero.  

\medskip
{\em Definition.\/} Let $\varphi_w:Q_w\prawo\eR$ be the linear form 
given by $\varphi_w(j)=1$ and $\varphi_w(y^\beta)=0$ for
$\beta\in\Delta-\{\bar\beta\}$. 
Let $\Phi_w:Q_w\times Q_w\prawo\eR$ be the corresponding symmetric 
bilinear form, $\Phi_w(f,g)=\varphi_w(fg)$.  
Let $M_w$ denote the matrix of $\Phi_w$ in the basis $y^\beta$, 
$\beta \in \Delta$.  
Let, as before, $S$ denote the multiplicative subset of $\Apis$ generated by
$g_{\beta_i}^{i}$.

\begin{lemma}  \label{cukier}
There is a symmetric matrix $M$ with entries in $S^{-1}\Apis$
such that $M_w=M(w)$ for $w\in W\setminus \Sigma$. $\Box$
\end{lemma}

  From now on we suppose that $F=(f_1\pkrop f_n):W\times\Rn\prawo\Rn$
is a polynomial mapping with $F(w;\gO)=\gO$ for every $w\in W$.
Denote
\[J=\frac{\partial(f_1\pkrop f_n)}{\partial(y_1\pkrop y_n)}
    \mbox{ and }J_w=J(w;\, \cdot\, ).\]

Let $\Ipis$ be the ideal in $\Apis[[y]]$ generated by $f_1\pkrop f_n$
and $I_w\subset\eR[[y]]$ that generated  by $f_1(w;\, \cdot\, )\pkrop
f_n(w;\, \cdot\, )$.  We assume that $\dim_R Q_w<\infty$ 
for every $w\in W$.  Hence, $\Delta$ and all $\Delta_w$ are finite. 

\begin{lemma} \label{kubek}
If $w\in W\setminus \Sigma$ then there is $0\neq\lambda_w\in\eR$ such that
$J_w=\lambda_w j$ in $Q_w$.
\end{lemma}

{\em Proof.\/}  By the Eisenbud-Levine Theorem \ref{el}, $J_w\neq 0$ 
in $Q_w$.  By Theorem \ref{mrowka}, 
$J_w=\sum_{\beta\in\Delta}\lambda_\beta y^\beta$
in $Q_w$.
Suppose, contrary to our claim, that $\lambda_{\beta'}\neq 0$ for a 
$\beta'<\bar\beta$.  Then, define a linear form 
$\psi:Q_w\prawo\eR$ by the formula $\psi(f)=f_{\beta'}\lambda_{\beta'}$, 
where $f=\sum_{\beta\in\Delta}f_\beta y^\beta\in Q_w$.  We show 
that the corresponding symmetric bilinear form $\Psi (f,g) = \psi (fg)$ 
is degenerate.  
For any $f\in Q_w$ we have
$\nu(fj)=\nu(f)+\nu(j)\geq\nu(j)=\bar\beta$.  Therefore, by Corollary 
\ref{rower}, $\psi(fj)=0$ for any $f\in Q_w$ and hence 
$\Psi(f,g)$ is degenerate.  On the other hand 
$\psi(J_w)=\lambda_{\beta'}^{2}>0$, and hence the existence of $\psi$ 
contradicts Theorem \ref{el}.  

Thus $\lambda_{\beta'}=0$ for every
$\beta' \in\Delta\setminus \{\bar\beta\}$ and we take $\lambda_w= 
\lambda_{\bar\beta}$. $\Box$ 
\medskip

In particular, by Theorem \ref{mrowka}, there is
$\lambda\in S^{-1}\Apis$ such that
$\lambda_w=\lambda(w)$ for $w\in W\setminus \Sigma$.\\[0.7em] 

{\em Definition.\/} Let 
$\psi_w=\lambda_w\varphi_w:Q_w\prawo\eR$.
Let $\Psi_w$ be the corresponding symmetric bilinear form.

\begin{proposition}  \label{trasa}
The forms $\psi_x$ and $\Psi_w$ defined above satisfy 
\begin{enumerate}
\item $\psi_w(J_w)>0$,
\item $\Psi_w$ is non-degenerate,
\item the entries and the determinant of the matrix of $\Psi_w$ in basis 
$y^\beta$, $\beta\in \Delta$, 
belong to $S^{-1}\Apis$. \end{enumerate}
\end{proposition} 

{\em Proof.\/} 
$\psi_w(J_w)=\lambda_w\varphi_w(J_w)=\lambda_{w}^{2}\varphi_w(j)=
  \lambda_{w}^{2}>0$, so the statement follows from the Eisenbud-Levine
Theorem \ref{el} and Lemma \ref{cukier}. $\Box$ \\ [0.4em]

Clearly multiplication by a positive scalar does not change the signature of
a symmetric matrix. So if we multiply the matrix of
$\Psi_w$ by the product of squares of the denominators of its entries 
we get

\begin{proposition}  \label{tara}
Assume that $W$ is irreducible.
Then there are a symmetric matrix $T$ with entries polynomials in $w\in W$ 
and a proper algebraic subset
$\Sigma\subset W$ such that for every $w\in W\setminus \Sigma$
\begin{enumerate}
\item $T(w)$ is non-degenerate,
\item $\signature \Psi_w=\signature \, T(w)$. $\Box$
\end{enumerate}
\end{proposition}

\bigskip
\section{Families of vector fields}\label{vectorfields}

\begin{lemma}  \label{marazm}
Let $F:W\times\Rn\prawo\Rn$ be a polynomial mapping. For any $w\in W$
let $F_w=F(w;\,\cdot\, ):\Rn\prawo\Rn$.  Suppose that 
for all $w\in W$, 
$\go\in\Rn$ is isolated in $F_{w}^{-1}(\go)$.   
(Hence $\deg_0 F_w$ is always well-defined.) 
Then there is a polynomial mapping $G:W\times\Rn\prawo\Rn$ such that
for every $w\in W$
\begin{enumerate}
\item $G_w:(\Rn,\gO)\prawo(\Rn,\gO)$ has an algebraically isolated
  zero at $\gO$,
\item $\deg_0F_w=\deg_0G_w$.
\end{enumerate}
\end{lemma}

{\em Proof.\/}  By the parametrized version of the {\L}ojasiewicz Inequality 
of \cite{fekak}, there
is $\alpha>0$ such that
\[\|  F_w(y)\| \geq C\|  y\| ^\alpha\]
for every $w\in W$ and $\| y\| <\delta$, where
$C=C(w)>0$ and $\delta=\delta(w)>0$ depend on $w$.

Choose an integer $k\gg 0$.  Define $G(w;y)=F(w;y)+(y_{1}^{k}\pkrop
y_{n}^{k})$. Let $G_{\bf C,w}:(\Cn,\gO)\prawo (\Cn,\gO)$ denote 
the complexification of $G_w$.  
Then, for every $w\in W$, $G_{\bf C,w}^{-1}(\gO)$ is a bounded complex 
algebraic set and hence finite. So $\gO$ is isolated in 
$G_{\bf C,w}^{-1}(\gO)$ and hence $G_w$ has an algebraically isolated zero 
at $\gO$.

We may assume that $k>\alpha$.  So if $w\in W$ and 
$y$ is close enough to the origin then
\[ \| tG_w(y)+(1-t)F_w(y)\| =
  \| F_w(y)+t(y_{1}^{k}\pkrop y_{n}^{k})\| \geq \]
\[ C\| y\|^{\alpha} - t \|(y_{1}^{k}\pkrop y_{n}^{k})\|  
 \geq \frac{C}{2} \|y\|^\alpha \]
for $0\leq t\leq 1$. Hence $\deg_0 F_w=\deg_0 G_w$ as required.  
$\Box$\\[0.7em]

\begin{lemma}\label{sloniczek} 
Under the assumptions of Lemma \ref{marazm}, if moreover $W$ is irreducible, 
then there exist a proper algebraic subset
$\Sigma\subset W$, an integer $\mu$, and polynomials
$q_1\pkrop q_t,q$ nowhere vanishing in $W\setminus \Sigma$ such that for every
$w\in W\setminus \Sigma$
\begin{enumerate}
\item $\deg_0 F_w={\sgn}\, q_1(w)\plkrop {\sgn}\, q_t(w)$,
\item $\deg_0 F_w\equiv \mu+1\pmod{2}$,
\item $\deg_0 F_w\equiv \mu+{\sgn}\,q(w)\pmod{4}$.
\end{enumerate}\end{lemma}

{\em Proof.\/} 
Let $F=(f_1\pkrop f_n)$, let $I_w$ denote the ideal in $\eR[[y]]$
generated by $f_1(w;\,\cdot\,),\ldots ,$ $f_n(w;\,\cdot\,)$ and let
$Q_w=\eR[[y]]/I_w$. By Lemma \ref{marazm} we may assume that each 
$F_w$ has an algebraically isolated zero at $\go$.  Let
\[J=\frac{\partial(f_1\pkrop f_n)}{\partial(y_1\pkrop y_n)}\]
and let $J_w$ denote the residue class of $J(w;\,\cdot\,)$ in $Q_w$.

Let $\psi_w:Q_w \to \eR$ be the linear form defined in section 3.  
By Proposition \ref{trasa}, $\psi_w$ satisfies the assumptions
of the Eisenbud-Levine Theorem \ref{el}. Hence the corresponding 
symmetric bilinear form $\Psi_w$ is non-degenerate
and $\deg_0 F_w=\mbox{signature}\,\Psi_w$.  
In particular, by Proposition \ref{tara}, there are a symmetric matrix 
$T$ with polynomial entries and a proper algebraic set $\Sigma'\subset W$ 
such that $T(w)$ is non-degenerate and
$\deg_0 F_w=$ signature$\,T(w)$ for every $w\in W\setminus \Sigma'$.  

Let $P_w(\lambda)=a_N\lambda^N+ a_{N-1}(w)\lambda^{N-1}\plkrop a_0(w)$, 
$a_N \equiv (-1)^N$,  
denote the characteristic polynomial of $T(w)$.  Clearly
its coefficients are polynomials in $w$ and $a_0(w)$ does not vanish
in $W\setminus \Sigma'$.
If $w\in W\setminus \Sigma'$ then all roots of $P_w$ are real and non-zero. 
Let$p_+(w)$ (resp. $p_-(w)$) denote the number of positive (resp. negative)
roots. Then
\[\mbox{signature}\,T(w)=p_+(w)-p_-(w), \]
and, by Lemma \ref{rosa}, it is easy to see
that there are a proper algebraic $\Sigma\subset W$, polynomials
$q_1\pkrop q_t,q$ nowhere vanishing on $W\setminus \Sigma$, and an integer 
$\mu$ such that 
\renewcommand{\labelenumi}{(\alph{enumi})}\begin{enumerate}
\item $\deg_0 F_w=p_+(w)-p_-(w)=\sign q_1(w)\plkrop \sign q_t(w)$,
\item $\deg_0 F_w\equiv \mu+\sign q(w)\pmod{4}$ \end{enumerate}
for every $w\in W\setminus \Sigma$ which completes the proof. $\Box$\\[1em]
\renewcommand{\labelenumi}{(\roman{enumi})}

Let $P$ be any non-negative polynomial with 
$P^{-1}(0)\cap W=\Sigma$.  Then
\[\sum \sign P(w)q_i(w)=\sum \sign q_i(w) \]
on $W\setminus \Sigma$ and 
\[\sum \sign P(w)q_i(w)=0\]  
on $\Sigma$.  Similarly, let $p_1\pkrop p_r$ be another set of polynomials.  
Then
\[\sum\sign p_j(w)+\sum\sign (-P(w)p_j(w))=0 \]
on $W\setminus \Sigma$ and 
\[\sum\sign p_j(w)+\sum\sign(-P(w)p_j(w))=\sum\sign p_j(w)\] 
on $\Sigma$.  Hence, by induction on $\dim W$ we get 

\begin{theorem}\label{slonik}
Let $W$ be a real algebraic set and let
$F:W\times\Rn\prawo\Rn$ be a polynomial mapping such that $\go$
is isolated in $F_{w}^{-1}(\go)$ for all $w\in W$. Then there are
polynomials $g_1\pkrop g_s$ such that for every $w\in W$
\[\deg_0 F_w={\sgn}\, g_1(w)\plkrop {\sgn}\, g_s(w). \hfil \Box\]
\end{theorem}

 % This is file: sferki.tex
 
\bigskip
\section{Families of algebraic sets}\label{families}
Let $X\subset W\times\Rn$ be a real algebraic set such that
$W\times\{\go\}\subset X$. There is a non-negative polynomial
$f:W\times\Rn\prawo\eR$ such that $X=f^{-1}(0)$. 
Denote $f_w(y)=f(w;y)$. Then $\go$ is contained in the set of critical points
of each $f_w.$
By the parametrized version of the {\L}ojasiewicz Inequality of \cite{fekak}, 
there is $\alpha>0$ such that for every $w\in W$ there are positive
$C=C(w)$ and $\delta=\delta(w)$ such that
\[ f_w(y)\geq C \|y\|^\alpha ,\]
for all critical points $y$ of $f_w$ with
$\|y\| <\delta$ and $f_w(y)\neq 0$.  

Let $k$ be an integer such that $2k>\alpha$.  Define  
\[g(w;y)=f(w;y)-\|y\|^{2k}\]
and let
\[G=\left( \frac{\partial g}{\partial y_1}\pkrop 
   \frac{\partial g}{\partial y_n}\right):W\times\Rn\prawo\Rn.\]
Clearly, $G$ is a polynomial family of vector fields
such that $G_w(\go)=\go$.

For every $w\in W$ let
$L(w)=\{y\in S_{r}^{n-1}\mid (w;y)\in X\}$, 
where $r>0$ is small. It is well-known that $L(w)$ is well-defined up 
to a homeomorphism.  Then $\chi(L(w))=1-\deg_0 G_w$.  
Indeed, this can be proven by an argument similar to that of proof of 
Lemma \ref{guzik}, if we replace $t^{2d+1}$ by $\|y\|^{2k}$, 
$P_\eta$ by the sphere $S_r$, and $\Sigma_\eta$ by the set of critical 
points of $f$ restricted to $S_r$, see \cite{szafraniec4} for the details.  
Therefore, Theorem \ref{slonik} implies

\begin{theorem} \label{antracyt} 
For all $w\in W$, 
$\Rn\ni\go$ is isolated in
$G^{-1}_{w}(\go)$ and
$\chi(L(w))=1-\deg_0 G_w$.   In particular, there are
polynomials $g_1\pkrop g_s$ such that for every $w\in W$
\[ \chi(L(w)) ={\sgn}\, g_1(w)\plkrop {\sgn}\, g_s(w). \quad \Box\] 
\end{theorem}

Similarly, let $S(w)=\{y\in S_{R}^{n-1}\mid (w;y)\in X\}$, 
where $R>0$ is very large.
$S(w)$ is well-defined up to a homeomorphism.

\begin{corollary} \label{rubin}
There is a polynomial family of vector fields
$H_w:\Rn\prawo\Rn$ such that $\Rn\ni\go$ is isolated in
$H_{w}^{-1}(\go)$ for all $w\in W$ and
$\chi(S(w))=1-\deg_0 H_w$. 
\end{corollary}

{\em Proof.\/} Let $d$ denote the degree of $f$, 
where as above, $f$ is a non-negative polynomial defining $X$.  
Then, there is a non-negative polynomial
$h:W\times\Rn\prawo\eR$ such that
$h(w;y)=\|y\|^{2d} f(w;y/\|y\|^2)$
for $y\neq \go$.  Clearly $h(w;\go)\equiv 0$ and $S(w)$ is
homeomorphic to
$L'(w)=\{y\in S_{r}^{n-1}\mid (w;y)\in h^{-1}(0)\},$ 
where $r>0$ is small.  
So the corollary follows from Theorem \ref{antracyt}. $\Box$

 % This is file: uzwart.tex

It is well-known (see, for instance, 
\cite{akbulutking4,benedettitognoli,bochnaketal})
that the single point
Aleksandrov compactification of a real algebraic set
is homeomorphic to a real algebraic set.  We shall
recall briefly the proof.

Suppose $X=\{y\in \Rn\mid f_1(y)=\ckrop=f_s(y)=0\}$,  where 
$f_1\pkrop f_s:\RndoR$ are polynomials of degree $\leq p-1$.  Set 
$h(y,y_{n+1})=y_{n+1}^{2}(f_{1}^{2}(y)\plkrop f_{s}^{2}(y))+(y_{n+1}-1)^{2}$, 
so that $h^{-1}(0)$ is homeomorphic to $X$ and $h$ is a non-negative
polynomial of degree $\leq 2p$.  
Put $y'=(y,y_{n+1})\in\Rn\times\eR$ and
$H(y')=\|  y'\| ^{4p}h(y'/\|  y'\| ^2)$.  
Then, it is easy to see that $H$ extends to a non-negative polynomial
on $\Rn\times\eR$ such that $H(\go,0)=0$ and
$H(y')=\|  y'\| ^{4p} +$ {\em monomials of lower degree\/}.  
Clearly $\tilde X=H^{-1}(0)$ is the single point compactification
of $X$ (If $X$ is compact then $\tilde X=X\amalg\{${\em point\/}$\}$).
Note that $t^{4p}H(y'/t)$ extends to a non-negative homogeneous
polynomial $f(y',t)$  on $\Rn\times\eR\times\eR$ of degree $4p$ such that
$f(y',0)=\|y'\|^{4p}$ and
$\tilde X$ is homeomorphic to
$\{y'\mid f(y',1)=0\}$.  Proceeding exactly in the same way we may prove 
the following parametrized
version of the above compactification method.

\begin{lemma} \label{rusznica}
Let $X\subset W\times\Rn$ be a real algebraic set. Then there
is a non-negative polynomial $f:W\times\eR^{n+1}\times\eR\prawo\eR$
such that for every $w\in W$ 
\begin{enumerate}
\item $f_w(y',t)=f(w;y',t):\eR^{n+1}\times\eR\prawo\eR$ is a non-negative
  homogeneous polynomial of degree $4p,$
\item $f_w(y',0)=\|  y'\| ^{4p},$
\item $\tilde X_w=\{y'\in\eR^{n+1}\mid f_w(y',1)=0\}$ is homeomorphic to
   the single point compactification of
   $X_w=\{y\in\Rn\mid (w;y)\in X\}$. $\Box$ 
\end{enumerate}
\end{lemma}

In particular, by Proposition \ref{guzik} we get

\begin{proposition} \label{arogant}
Let $X\subset W\times\Rn$ be a real algebraic set.
Then there is a polynomial family of vector fields
$F_w:\RndoRn$ such that for every $w\in W$ 
\begin{enumerate}
\item $F_w(\go)=\go,$
\item $\go$ is isolated in $F_{w}^{-1}(\go),$
\item $\chi(\tilde X_w)=1-\deg_0 F_w .\ \Box$ \end{enumerate}
\end{proposition}

Let $S(w)=X_w\cap S_{R}^{n-1},$ where $R>0$ is sufficiently large. 
Then it is easy
to check that
\[\chi(X_w)=\chi(\tilde X_w)+\chi (S(w))-1.\]
By \ref{arogant}, \ref{rubin}, \ref{slonik}, and since 
$\sign a+\sign b\equiv \sign (ab)+1\pmod{4}$, provided $a\neq 0$ and 
$b\neq 0$, we get

\begin{theorem} \label{wtyczka} 
Let $X\subset W\times \Rn$ be a real algebraic set.  Then 
there are polynomials $g_1\pkrop g_s$ on $W$ such that
\[\chi(X_w)={\sgn}\,g_1(w)\plkrop {\sgn}\,g_s(w).\]
In particular, if $W$ is irreducible, then there are a proper algebraic
subset $\Sigma\subset W$, an integer $\mu$, and a polynomial
$g$ nowhere vanishing in $W-\Sigma$ such that for every
$w\in W-\Sigma$ 
\[\chi(X_w)\equiv \mu+ {\sgn}\,g(w)\pmod{4},\]
In particular $\chi(X_w)\equiv\mu+1\pmod{2}$. $\Box$ 
\end{theorem}

\bigskip
\section{Algebraically constructible functions}\label{constructible}

Let $W$ be a real algebraic set.  An integer-valued function 
$\varphi:W\to \mathbf Z$ is called ({\it semialgebraically}) 
{\it constructible}  
if it admits a presentation as a finite sum
\begin{equation}\label{constr}
\varphi = \sum m_i  \one_{W_i}, 
\end{equation}
where for each $i$, $W_i$ is a semialgebraic subset of $W$,
$\one_{W_i}$ is the characteristic function of $W_i$, and $m_i$ is an
integer.  Constructible functions, well-known in complex domain, were studied 
in real algebraic set-up by Viro \cite{viro}, and in sub-analytic set-up 
by Kashiwara and Schapira \cite{kashiwaraschapira, schapira}.  
If the support  of constructible function $\varphi$
is compact, then we may choose  all $W_i$ in (\ref{constr}) compact.  
Then, cf. \cite{viro, schapira, mccroryparusinski2}, 
the {\it Euler integral} of $\varphi$ is defined as  
\[\int \varphi = \sum m_i \chi (W_i) . \]
It follows from the additivity of Euler characteristic that
the Euler integral is well-defined and does not depend on the presentation 
(\ref{constr}) of $\varphi$, provided all $W_i$ are compact.  
Let $f:W\to Y$ be a (continuous) semialgebraic map of real algebraic
sets, $\varphi$ a constructible function on $W$ and suppose that
 $f:W\to Y$ restricted to the support of $\varphi$ is proper.  
 Then the {\it direct image} $f_*\varphi$ is
 given by the  formula
\[f_*\varphi (y) = \int_{f^{-1}(y)} \varphi ,\]
where by $\int_{f^{-1}(y)} \varphi$ we understand the Euler integral 
of $\varphi$ restricted to $f^{-1}(y)$.  It follows from the existence 
of a stratification of $f$ that $f_*\varphi$ is a constructible function
on $Y$.  

Another more restrictive class of constructible functions, was introduced 
in \cite{mccroryparusinski2} in order to study local topological properties 
of real algebraic sets.    An integer-valued function 
$\varphi:W \to \mathbf Z$ is called {\it algebraically
constructible} if there exists a finite collection of algebraic sets
$Z_i$, regular proper morphisms ${f_i}:Z_i \to W$, and integers $m_i$, 
 such that 
\[ \varphi = \sum m_i {f_i}_* \one_{Z_i} \]

It is obvious that every algebraically constructible function is 
semialgebraically constructible but the converse is false for 
$\dim W >0$.  For instance, a constructible function on $\mathbf R$ 
is algebraically constructible if and only if it is is generically 
constant mod 2.  The reader may consult \cite{mccroryparusinski2} for 
other examples.  As a consequence of section \ref{families} we obtain 
the following simple decription of algebraically constructible functions.  

\begin{theorem}\label{key}
Let $W$ be a real algebraic set.  Then $\varphi:W\to \mathbf Z$ 
is algebraically constructible if and only if there exist polynomial 
functions $g_1\pkrop g_s$ on $W$ such that
\[\varphi(w)=\sgn g_1(w)\plkrop \sgn g_s(w).\]
\end{theorem}

{\em Proof.\/}  It is easy to see that the sign of a polynomial function
$g$ on $W$ defines 
an algebraically constructible function.  
Indeed, let 
$\widetilde W = \{(w, t) \in W\times \eR \ |\ g(w) = t^2\}$ and   
let $\pi : \widetilde W\to W$ denote the standard projection.  
Then $\sgn f = \pi_* \one_{\widetilde W} - \one_W$ is algebraically 
constructible.  

The opposite implication follows from Theorem \ref{wtyczka}. $\Box$ \par

\begin{corollary} \label{gniazdko} \par 
\begin{enumerate} 
\item 
Let $F:W\times\Rn\prawo\Rn$ be a polynomial mapping satisfying the 
assumptions of \ref{slonik}.  
Then $w\prawo \deg_0 F_w$ is an algebraically constructible function on $W$.
\item 
Let $X_w$ be an algebraic family of affine real algebraic sets parametrized 
by $w\in W$ as in \ref{wtyczka}.  Then $w\prawo \chi(X_w)$ is an algebraically 
constructible function on $W$. $\Box$ 
\end{enumerate}
\end{corollary}

The next corollary is virtually equivalent to the main result of 
\cite{costekurdyka2}.  

\begin{corollary}\label{discriminant} 
Let $\varphi$ be an algebraically constructible function on 
an irreducible real algebraic set $W$.  
Then there exist a proper real algebraic subset 
$\Sigma\subset W$, an integer $\mu$, and a polynomial
$g$ on $W$, such that $g$ does not vanish on $W\setminus \Sigma$ and  
\[ \varphi(w)\equiv \mu+{\sgn}\,g(w) \pmod{4}\] 
for $w\in W-\Sigma$.
In particular, for such $w$, $\varphi(w)\equiv\mu+1\pmod{2}$.  
\end{corollary} 

{\em Proof.\/} Let $g_1\pkrop g_s$ be polynomials given by \ref{key}.  
We may suppose that all of them are not identically equal to zero.  
Since $\sign a+\sign b\equiv \sign (ab)+1\pmod{4}$, for $a$ and $b$ non-zero, 
the polynomial  
$g= g_1\cdots g_s$ satisfies the statement.  
This ends the proof. $\Box$ \par

Let $\varphi$ be a constructible function on $W$.  Following 
\cite {mccroryparusinski2}  we define
the {\it link} of $\varphi$ as the constructible function on $W$ given by
\[\Lambda \varphi (w) = \int_{S(w,\varepsilon)} \varphi , \]
where
$\varepsilon>0$ is sufficiently small, and $S(w, \varepsilon)$ denotes the
$\varepsilon$-sphere centered at $w$.  It is easy to see that $\Lambda\varphi$ 
is well defined and independent of the embedding of $W$ in $\bf R^n$.  
Then the duality operator $D$ on constructible functions, introduced by 
Kashiwara and Schapira in \cite {kashiwaraschapira, schapira},  
satisfies
\[ D \varphi = \varphi - \Lambda \varphi.  \]
As shown in \cite{mccroryparusinski2} the following general statement 
generalizes various previously known restrictions on local topological 
properties of real algebraic sets.  In particular it implies Akbulut 
and King's numerical conditions of \cite{akbulutking} and the 
conditions modulo 4, 8, and 16 of Coste and Kurdyka 
\cite {coste2, costekurdyka1} generalized in \cite{mccroryparusinski1}.  

\begin{theorem}\label{link} 
Let $\varphi$ be an algebraically constructible function on a real 
algebraic set $W$.  Then $\frac 1 2 \Lambda \varphi$ is 
integer-valued and algebraically constructible. 
\end{theorem} 

The above theorem was proven in \cite {mccroryparusinski2} using the 
resolution of singularites.  As we show below it is a simple conseqence of 
Theorem \ref{key}.  

{\em Proof.\/} 
We begin the proof by some preparatory observations.    

\begin{lemma}\label{limits}
$W$ be a real algebraic set and let $\gamma$ be an algebraically 
constructible function on $W\times \mathbf R$.  Then
\[ \psi_+ (w)= \lim_{t\to 0_+} \gamma(w,t),  \, \, 
\psi_- (w)=  \lim_{t\to 0_+} \gamma(w,-t),  \, \,
\psi (w)=  \frac 1 2 (\psi_+ (w) - \psi_- (w)) \]
are integer-valued and algebraically constructible  on 
$W_0 = W\times \{0\}$.  
\end{lemma}

{\em Proof.\/} 
We show the lemma for $\psi$.  The proofs for $\psi_+$ and $\psi_-$ are  
similar.  

We proceed  by induction on $\dim W$.  Without loss of generality we 
may assume that that $W$ is affine and irreducible.  
We shall show that the statement of lemma holds generically on $W_0$, 
that is to say there exists a proper algebraic subset $W'_0$ of $W_0$ and an 
algebraically constructible function $\psi'$ on $W_0$ which 
equals $\psi$ in the complement of $W'_0$.  Then the lemma follows 
 from the inductive assumption since $\dim W'_0 < \dim W_0$.  

By Theorem \ref{key} we may assume that $\gamma = \sign g$, where $g(w,t)$ 
is a polynomial function on $W\times \mathbf R$.  We may also assume that 
$g$ does not vanish identically, and then there exists 
a nonnegative integer $k$ such that 
\[ g(w,t) = t^k h(w,t), \]
where $h(w,t)$ is a polynomial function on $W\times \mathbf R$ not vanishing 
identically on $W\times \{0\}$.  Then, in the complement of 
$W'_0 = \{w| h(w,0)=0\}$, 
either $\psi(w) = \sgn h(w,0)$ for $k$ odd or $\psi(w) =0$ for $k$ even 
satisfies the statement.  
This ends the proof of lemma. $\Box$  

Let $\widetilde W = \{(w,y,t)\in W\times W\times\mathbf R|\, \|w-y\|^2=t \}$ 
and let $\pi :\widetilde W\to W\times \mathbf R$ be given by  
$\pi (w,y,t) = (w,t)$.  Let $\tilde \varphi (w,y,t)= \varphi (y)$.  
Then $\tilde \varphi$ is algebraically 
constructible and hence  $\gamma = \pi_* \tilde \varphi$ is an 
algebraically constructible function on $W\times \mathbf R$ and 
\[ \lim_{t\to 0_+} \gamma(w,t) =  \Lambda \varphi (w) . \]
Since $\gamma (w,t)=0$ for $t<0$   
\[ \frac 1 2 \Lambda \varphi (w) = \frac 1 2 \lim_{t\to 0_+} 
(\gamma(w,t)-\gamma(w,-t)) 
\]
is algebraically constructible by Lemma \ref{limits}.  
This ends the proof of Theorem \ref{link}. $\Box$ 
\medskip

Suppose that  $f:W\to \eR$ is regular and let  $w\in W_0 = f^{-1} (0)$.  
Then we define the {\it positive}, resp.~{\it negative}, {\it Milnor fibre}
of $f$ at $w$ by
\[   F_f^+(w)  =  B(w, \varepsilon )\cap f^{-1} (\delta) \]   
\[ F_f^-(w) =  B(w, \varepsilon )\cap f^{-1} (-\delta) ,\]
where  $B(w, \varepsilon )$ is the ball of radius
$\varepsilon$ centered at $w$ and $0<\delta\ll \varepsilon\ll 1$.
Let $\varphi$ be an algebraically constructible function on $W$.  
Following \cite{mccroryparusinski2} we define the {\it positive} (resp.~{\it
negative}) {\it specialization} of $\varphi$ with respect to $f$ by
\[ (\Psi^+_f \varphi) (w) = \int_{F_f^+(w)} \varphi, \quad   (\Psi^-_f
\varphi) (w) = \int_{F_f^-(w)} \varphi . \]
It is easy to see that both specializations are well-defined and that they 
are constructible functions supported in $W_0$.  Moreover, as shown in 
\cite{mccroryparusinski2},
they are also algebraically costructible.  
We present below an alternative 
proof of this fact.  

\begin{theorem}\label{specialization} 
Let  $f:W\to \eR$ be a regular function on a real algebraic set 
$W$.  Let $\varphi$ be an algebraically constructible function on $W$.  Then 
$\Psi^+_f \varphi$, $\Psi^-_f \varphi$, and  
$\frac 1 2 (\Psi^+_f \varphi - \Psi^-_f \varphi)$ are integer valued and 
algebraically constructible.  
\end{theorem} 

{\em Proof.\/} 
The proof is similar to that of Theorem \ref{link}.  Since the Milnor fibres 
are defined not only by equations but also by inequalities we use the 
following auxiliary construction.  

Let $\widetilde W = \{(w,y,t,r,s)\in W\times W\times\mathbf R^3|\, 
\|w-y\|^2 + t^2=r, f(y)=s \}$ 
and let $\pi :\widetilde W\to W\times \eR ^2$ be given by  
$\pi (w,y,t,r,s) = (w,r,s)$.  Note that for $w\in W_0$, $0<s\ll r\ll 1$, 
$\tilde F = \pi^{-1} (w,r,s)$ is a double cover of the Milnor fibre 
$F=F^+_f(w)$ branched along its boundary 
$\partial F = S(w,\sqrt r)\cap f^{-1}(s)$.  
Hence $\chi (\tilde F) = 2\chi (F) - \chi (\partial F)$. 
Let $\tilde \varphi (w,y,t,r,s)= \varphi (y)$.  Then 
\[ \Psi_f^+ \varphi (w) = \frac 1 2 (\int_{\tilde F}\tilde \varphi + 
\int_{\partial F} \varphi) = \frac 1 2 \pi_* 
(\tilde \varphi + \tilde \varphi|_{t=0}) (w,r,s), \]  
for $0<s\ll r\ll 1$.  Clearly an analogous formula holds for 
$\Psi_f^- \varphi (x)$.

Let  $\gamma = \pi_* (\tilde \varphi + \tilde \varphi|_{t=0})$.  
Then $\gamma (w,r,s)$ is algebraically constructible and $\gamma(w,r,s)=0$ 
for $r<0$.  Hence, by Lemma \ref{limits}, the following functions are 
algebraically  constructible 
\[ \Psi_f^{\pm} \varphi = \frac 1 2 \lim_{r\to 0_+} \lim_{s\to 0_+} 
\gamma(w,r,\pm s) , \]
\[ \frac 1 2 (\Psi_f^{+} - \Psi_f^-) \varphi = 
\frac 1 4 \lim_{r\to 0_+} \lim_{s\to 0_+} 
(\gamma(w,r,s)-\gamma(w,r,-s)), \]
as required. $\Box$

\end{document}